\documentclass[12pt,reqno]{amsart}

\usepackage{amssymb,amsmath}
\usepackage{amsfonts}
\usepackage{amsthm}
\usepackage{mathrsfs}
\usepackage{bm}
\usepackage{graphicx}

\theoremstyle{plain} 
\newtheorem{theorem}{Theorem}
\newtheorem{corollary}[theorem]{Corollary}

\newtheorem{lemma}[theorem]{Lemma}

\theoremstyle{definition} 
\newtheorem{definition}[theorem]{Definition}
\newtheorem{remark}[theorem]{Remark}
\newtheorem{example}[theorem]{Example}
\newcommand{\R}{\ensuremath{\mathbb{R}}}

\newcommand{\N}{\ensuremath{\mathbb{N}}}
\newcommand{\C}{\ensuremath{\mathbb{C}}}

\numberwithin{equation}{section}
\numberwithin{theorem}{section}

\small\normalsize

\begin{document}

\title[Best Ulam constants for damped linear oscillators]{Best Ulam constants for damped linear oscillators with variable coefficients}
\author[Anderson]{Douglas R. Anderson} 
\address{Department of Mathematics, 
         Concordia College, 
         Moorhead, MN 56562 USA}
\email{andersod@cord.edu}
\author[Onitsuka]{Masakazu Onitsuka}
\address{Department of Applied Mathematics, 
         Okayama University of Science, 
         Okayama, 700-0005, Japan}
\email{onitsuka@ous.ac.jp}
\author[O'Regan]{Donal O'Regan}
\address{School of Mathematical and Statistical Sciences, 
         University of Galway, 
         Galway, Ireland}
\email{donal.oregan@nuigalway.ie}

\keywords{Ulam stability; best Ulam constant; damped linear oscillator; second-order linear differential equation; Lane-Emden differential equation; Riccati equation; variable coefficient.}
\subjclass[2020]{34D10, 34D20, 34A30}

\begin{abstract}
This study uses an associated Riccati equation to study the Ulam stability of non-autonomous linear differential vector equations that model the damped linear oscillator. In particular, the best (minimal) Ulam constants for these non-autonomous linear differential vector equations are derived. These robust results apply to vector equations with solutions that blow up in finite time, as well as to vector equations with solutions that exist globally on $(-\infty,\infty)$. Illustrative, non-trivial examples are presented, highlighting the main results. 
\end{abstract}

\maketitle\thispagestyle{empty}


\section{Introduction}

Let $I$ be an interval of $\R$, and let $n\in \N$. In this paper, we consider the second-order non-autonomous linear differential vector equation
\begin{equation}
 \alpha(t)\bm{x}'' + \beta(t)\bm{x}' + \gamma(t)\bm{x} = \bm{f}(t),
 \label{oscillator1}
\end{equation}
where $\alpha$, $\beta$, $\gamma: I \to \C$ are continuous scalar functions, and $\bm{f}: I \to \C^n$ is a continuous vector function, and $\bm{x}\in C^2(I,\C^n)$. Moreover we assume that $\alpha(t)\ne 0$ for all $t \in I$. As is well known, when the coefficients $\alpha$, $\beta$, and $\gamma$ of \eqref{oscillator1} are all positive constants, it is called a forced, damped linear oscillator, and has various applications in for example mechanical engineering. In many mathematical models, constant coefficients are determined through statistical processing based on actual phenomena. However, it is also true that constants may not always be sufficient due to factors such as various changes in temperature, humidity, etc. due to seasonal changes. In addition, there are phenomena for which variable coefficients are appropriate. For this reason, it can be expected that considering the non-autonomous equation \eqref{oscillator1} will be useful in applications. So what is the biggest problem when considering mathematical models? It goes without saying that much of it has to do with perturbation. There will always be some kind of error between the actual phenomenon and the mathematical model, and depending on the situation, it may lead to a large error as time passes. Therefore, ensuring that the error does not grow large is a very important application issue. Concepts such as ``Ulam stability'' originating from the field of functional equations (see, \cite{AgaXuZha,BrzPopRasXu,BrzPopRasXu2,Cieplinski1,Cieplinski2}) and ``shadowing'' originating from the field of dynamical systems (see, \cite{BacDra1,BacDra2,BacDraOniPit,BacDraPitSin}) are known to tackle such problems. Both of these are very similar, and this study will use the former definition to proceed with the discussion. 


\begin{definition}\label{Ulamdef}
Equation \eqref{oscillator1} is Ulam stable on $I$, if there exists a constant $K>0$ such that for every $\varepsilon>0$, and for every twice continuously differentiable function $$\bm{\xi} \in C^2(I, \C^n)$$ satisfying
\[ \sup_{t \in I}\|\alpha(t)\bm{\xi}'' + \beta(t)\bm{\xi}' + \gamma(t)\bm{\xi} - \bm{f}(t)\| \le \varepsilon, \]
there exists a solution $\bm{x} \in C^2(I, \C^n)$ of \eqref{oscillator1} such that
\[ \sup_{t \in I}\|\bm{\xi}(t)-\bm{x}(t)\| \le L \varepsilon. \]
We call such an $L$ an Ulam constant for \eqref{oscillator1} on $I$. Here $\|\bm{y}\|$ represents the norm of $\bm{y} \in \C^n$. 
\end{definition}

Since the concept of Ulam stability was introduced into differential equations much research has been done. There are many results even limited to second-order differential equations, and they have been actively studied in recent years. For example, see \cite{Dragicevic1,GavJunLi,LiHua,LiShe,WenWanO'Re,YanMen,ZadAlaXuDon}. Focusing on second-order differential equations with constant coefficients, strong results have been obtained. In 2020, Baias and Popa \cite{BaiPop1} investigated the details of the solutions of the constant coefficients equation
\begin{equation}
 x'' + a_1x' + a_2x = 0,
 \label{oscillator2}
\end{equation}
and its perturbed equations, and derived the minimum Ulam constant, where $a_1$ and $a_2$ are complex-valued constants, and $(X, \| \cdot\|)$ is a Banach space over the field $\C$, and $x \in C^2(\R, X)$. They analyzed the Ulam stability based on case division according to the signs of the real parts of the roots of the characteristic equations corresponding to \eqref{oscillator2}. Finding the minimum Ulam constant is a very important problem for applications. If we recall the words mathematical model and phenomenon from earlier, the minimum Ulam constant means the minimum error between the mathematical model and phenomenon. Establishing Ulam stability remains an important qualitative property, but finding the minimum Ulam constant quantitatively specifies the minimum error between the approximate and true solutions. The minimum Ulam constant is usually called the best Ulam constant, and research on the best Ulam constant for differential equations, difference equations, and functional equations with constant coefficients has progressed in recent years (see, \cite{AndOni1,AndOni2,AndOni3,BaiPop2,BaiPop3,PopRas1}). On the other hand, for equations with variable coefficients, there are few studies establishing even the Ulam stability, let alone the best Ulam constants. As far as the authors are aware, the best Ulam constants have been derived only for periodic equations (see, \cite{AndOniRas,FukuOni1,FukuOni2}). The purpose of this paper is to derive the minimum (best) Ulam constants for the non-autonomous equation \eqref{oscillator1}, which is not limited to periodic coefficients. 

In 2020, C\u{a}dariu, Popa and Ra\c{s}a \cite{CadPopRas} investigated the Ulam stability of the second-order non-autonomous linear differential scalar equation
\begin{equation}
 x'' + \beta(t)x' + \gamma(t)x = 0,
 \label{oscillator3}
\end{equation}
where $\beta \in C^1(I, \R)$, $\gamma \in C(I, \R)$ are real-valued scalar functions. They used the existence of solutions to the initial value problem of a particular Riccati equation to give a result that guarantees the Ulam stability of \eqref{oscillator3}. Later, this result was extended to the case where \eqref{oscillator1} was limited to scalar equations and real-valued coefficients (see, \cite{Onitsuka2}). However, unfortunately, the best Ulam constants were not obtained in the above two results. This study attempts to analyze \eqref{oscillator1} using their idea that the existence of solutions to a Riccati equation is useful in analyzing the Ulam stability of second-order linear differential equations. It should be noted here that the present study does not extend their results. This work does not require continuous differentiability for the coefficient $\beta$. It is enough to assume continuity. In addition, the Riccati equation used in this study is also different from the one proposed by them, and the method is completely different. Therefore, the statements of the obtained theorems are also different from theirs. In this study, by proposing a new method, we succeeded in deriving the best Ulam constants for \eqref{oscillator1}. 

This paper is organized as follows. In Section 2, we show that when we guarantee the existence of a solution to a certain Riccati equation, we can use it to describe the solution to the initial value problem of \eqref{oscillator1}. In Section 3, we give the main theorem and its proof. Section 4 gives the result of deriving the best Ulam constants, which is the goal of this study. In Section 5, we present various non-trivial examples centering on Lane-Emden differential equations.


\section{Representation of solution}

In this section, we show that, given the existence of a solution to a certain Riccati equation, we can use it to express the general solution of \eqref{oscillator1}. 

\begin{lemma}\label{genesol}
Suppose that $\alpha(t)\ne 0$ for all $t \in I$, and there exists a solution $\rho\in C^1(I, \C)$ of the Riccati equation
\begin{equation}
 \alpha(t)(\rho'+\rho^2) + \beta(t)\rho + \gamma(t) = 0.
 \label{Riccati}
\end{equation}
Then the solution of \eqref{oscillator1} with $\bm{x}(t_0)=\bm{x}_0$ and $\bm{x}'(t_0)=\bm{x}'_0$ is given by
\[ \bm{x}(t) = \left[ \bm{x}_0+\int_{t_0}^{t} \left(\bm{x}'_0-\rho(t_0)\bm{x}_0+\int_{t_0}^{s}\frac{e^{\int_{t_0}^{\mu}\left(\rho(\nu)+\frac{\beta(\nu)}{\alpha(\nu)}\right)d\nu}}{\alpha(\mu)}\bm{f}(\mu)d\mu \right) e^{-\int_{t_0}^{s}\left(2\rho(\mu)+\frac{\beta(\mu)}{\alpha(\mu)}\right)d\mu}ds \right] e^{\int_{t_0}^{t}\rho(s)ds} \]
for $t \in I$, where $t_0 \in I$. 
\end{lemma}

\begin{proof}
Assume that $\alpha(t)\ne 0$ for all $t \in I$. Let $\rho(t)$ be a solution of \eqref{Riccati} on $I$. Then we see that the function
\[ \bm{y}(t) = \bm{c}e^{\int_{t_0}^{t}\rho(s)ds}, \quad \bm{c} \in \C^n \]
is a solution to the damped linear oscillator
\[ \alpha(t)\bm{y}'' + \beta(t)\bm{y}' + \gamma(t)\bm{y} = \bm{0} \]
on $I$. Now we use the reduction of order method. Letting
\[ \bm{x}(t) = \bm{z}(t)e^{\int_{t_0}^{t}\rho(s)ds}, \]
we have
\[ \bm{x}'(t) = \left(\bm{z}'(t)+\rho(t)\bm{z}(t)\right)e^{\int_{t_0}^{t}\rho(s)ds}, \]
and
\[ \bm{x}''(t) = \left[\bm{z}''(t)+2\rho(t)\bm{z}'(t)+\left(\rho'(t)+\rho^2(t)\right)\bm{z}(t)\right]e^{\int_{t_0}^{t}\rho(s)ds}. \]
Substituting these into \eqref{oscillator1} and using \eqref{Riccati}, we obtain
\[ \alpha(t)\bm{z}''(t)+\left(2\alpha(t)\rho(t)+\beta(t)\right)\bm{z}'(t) = \bm{f}(t)e^{-\int_{t_0}^{t}\rho(s)ds}. \]
Since $\alpha(t)\ne 0$ for all $t \in I$, we have
\[ \left(\bm{z}'(t)e^{\int_{t_0}^{t}\left(2\rho(s)+\frac{\beta(s)}{\alpha(s)}\right)ds}\right)' = \frac{e^{\int_{t_0}^{t}\left(\rho(s)+\frac{\beta(s)}{\alpha(s)}\right)ds}}{\alpha(t)}\bm{f}(t). \]
This implies that
\[ \bm{z}'(t) = \left(\bm{z}'(t_0)+\int_{t_0}^{t}\frac{e^{\int_{t_0}^{s}\left(\rho(\mu)+\frac{\beta(\mu)}{\alpha(\mu)}\right)d\mu}}{\alpha(s)}\bm{f}(s)ds \right)e^{-\int_{t_0}^{t}\left(2\rho(s)+\frac{\beta(s)}{\alpha(s)}\right)ds}, \]
and that
\[ \bm{z}(t) = \bm{z}(t_0)+\int_{t_0}^{t}\left(\bm{z}'(t_0)+\int_{t_0}^{s}\frac{e^{\int_{t_0}^{\mu}\left(\rho(\nu)+\frac{\beta(\nu)}{\alpha(\nu)}\right)d\nu}}{\alpha(\mu)}\bm{f}(\mu)d\mu \right)e^{-\int_{t_0}^{s}\left(2\rho(\mu)+\frac{\beta(\mu)}{\alpha(\mu)}\right)d\mu}ds. \]
From
\[ \bm{z}(t_0) = \bm{x}(t_0) = \bm{x}_0, \quad \bm{z}'(t_0) = \bm{x}'(t_0)-\rho(t_0)\bm{x}(t_0) = \bm{x}'_0-\rho(t_0)\bm{x}_0, \]
we obtain
\[ \bm{x}(t) = \left[ \bm{x}_0+\int_{t_0}^{t} \left(\bm{x}'_0-\rho(t_0)\bm{x}_0+\int_{t_0}^{s}\frac{e^{\int_{t_0}^{\mu}\left(\rho(\nu)+\frac{\beta(\nu)}{\alpha(\nu)}\right)d\nu}}{\alpha(\mu)}\bm{f}(\mu)d\mu \right) e^{-\int_{t_0}^{s}\left(2\rho(\mu)+\frac{\beta(\mu)}{\alpha(\mu)}\right)d\mu}ds \right] e^{\int_{t_0}^{t}\rho(s)ds} \]
for $t \in I$. This completes the proof.
\end{proof}

\begin{remark}
If we assume $\alpha(t)\ne 0$ almost everywhere $t \in I$ and with all the obvious functions in the integrands in the expression two lines after \eqref{Riccati} in $L^1$, then we can get a solution of \eqref{oscillator1} in $W^{1,1}(I, \C^n)$. 
\end{remark}


\section{Ulam stability}

The first theorem of this paper is as follows. 


\begin{theorem}\label{UlamS}
Let $I$ be either $(\tau,\sigma)$, $(\tau,\sigma]$, $[\tau,\sigma)$ or $[\tau,\sigma]$, where $-\infty \le \tau < \sigma \le \infty$. Suppose that $\alpha(t)\ne 0$ for all $t \in I$, and there exists a solution $\rho: I \to \C$ of \eqref{Riccati}. Let $\Re(z)$ be the real part of $z\in\C$. Then the following (i), (ii) and (iii) below hold:
\begin{itemize}
  \item[(i)] if the functions
\begin{equation}
 f_{1}(t):= \int_{t}^{\sigma} \frac{e^{\int_t^s \Re\left(\rho(\mu)+\frac{\beta(\mu)}{\alpha(\mu)}\right) d\mu}}{|\alpha(s)|}ds
 \label{function1}
\end{equation}
and
\begin{equation}
 f_{2}(t):= \int_{t}^{\sigma} e^{-\int_t^s \Re(\rho(\mu)) d\mu}ds
 \label{function2}
\end{equation}
exist for all $t \in I$, and $\sup_{t \in I}f_{1}(t) < \infty$ and $\sup_{t \in I}f_{2}(t) < \infty$ hold. Then \eqref{oscillator1} is Ulam stable on $I$, with an Ulam constant
\[ L_{1} := \sup_{t \in I} \int_{t}^{\sigma} \left(\int_{s}^{\sigma} \frac{e^{\int_s^\mu \Re\left(\rho(\nu)+\frac{\beta(\nu)}{\alpha(\nu)}\right) d\nu}}{|\alpha(\mu)|}d\mu\right) e^{-\int_t^s \Re(\rho(\mu)) d\mu}ds; \]
  \item[(ii)] if the functions $f_{1}(t)$ and
\begin{equation}
 f_{3}(t):= \int_{\tau}^{t} e^{\int_s^t \Re(\rho(\mu)) d\mu}ds
 \label{function3}
\end{equation}
exist for all $t \in I$, and $\sup_{t \in I}f_{1}(t) < \infty$ and $\sup_{t \in I}f_{3}(t) < \infty$ hold, where $f_{1}(t)$ is given by \eqref{function1}. Then \eqref{oscillator1} is Ulam stable on $I$, with an Ulam constant
\[ L_{2} := \sup_{t \in I} \int_{\tau}^{t} \left(\int_{s}^{\sigma} \frac{e^{\int_s^\mu \Re\left(\rho(\nu)+\frac{\beta(\nu)}{\alpha(\nu)}\right) d\nu}}{|\alpha(\mu)|}d\mu\right) e^{\int_s^t \Re(\rho(\mu)) d\mu}ds; \]
  \item[(iii)] if the functions $f_{3}(t)$ and
\begin{equation}
 f_{4}(t):= \int_{\tau}^{t} \frac{e^{-\int_s^t \Re\left(\rho(\mu)+\frac{\beta(\mu)}{\alpha(\mu)}\right) d\mu}}{|\alpha(s)|}ds
 \label{function4}
\end{equation}
exist for all $t \in I$, and $\sup_{t \in I}f_{3}(t) < \infty$ and $\sup_{t \in I}f_{4}(t) < \infty$ hold, where $f_{3}(t)$ is given by \eqref{function3}. Then \eqref{oscillator1} is Ulam stable on $I$, with an Ulam constant
\[ L_{3} := \sup_{t \in I} \int_{\tau}^{t} \left(\int_{\tau}^{s} \frac{e^{-\int_\mu^s \Re\left(\rho(\nu)+\frac{\beta(\nu)}{\alpha(\nu)}\right) d\nu}}{|\alpha(\mu)|}d\mu\right) e^{\int_s^t \Re(\rho(\mu)) d\mu}ds. \]
\end{itemize}
\end{theorem}

\begin{proof}
Assume that $\alpha(t)\ne 0$ for all $t \in I$. Assume also that there exists a solution $\rho: I \to \C$ of \eqref{Riccati}. Let $\varepsilon>0$ be given, and let the twice continuously differentiable function $\bm{\xi}: I \to \C^n$ satisfy
\[ \sup_{t \in I}\|\alpha(t)\bm{\xi}'' + \beta(t)\bm{\xi}' + \gamma(t)\bm{\xi} - \bm{f}(t)\| \le \varepsilon. \]
Define
\[ \bm{g}(t) := \alpha(t)\bm{\xi}'' + \beta(t)\bm{\xi}' + \gamma(t)\bm{\xi} - \bm{f}(t) \]
for $t \in I$. Then we have $\sup_{t \in I}\|\bm{g}(t)\| \le \varepsilon$. Let $\bm{p}(t)$ be a solution to \eqref{oscillator1} on $I$, and let $\bm{q}(t):=\bm{\xi}(t)-\bm{p}(t)$ for $t \in I$. Then $\bm{q}(t)$ is a solution to the equation
\[ \alpha(t)\bm{q}'' + \beta(t)\bm{q}' + \gamma(t)\bm{q} = \bm{g}(t) \]
for $t \in I$. Therefore, by Lemma \ref{genesol}, we see that the function $\bm{q}(t)$ is expressed as
\begin{equation}
 \bm{q}(t) = \left[ \bm{q}_0+\int_{t_0}^{t} \left(\bm{q}'_0-\rho(t_0)\bm{q}_0+\int_{t_0}^{s}\frac{e^{\int_{t_0}^{\mu}\left(\rho(\nu)+\frac{\beta(\nu)}{\alpha(\nu)}\right)d\nu}}{\alpha(\mu)}\bm{g}(\mu)d\mu \right) e^{-\int_{t_0}^{s}\left(2\rho(\mu)+\frac{\beta(\mu)}{\alpha(\mu)}\right)d\mu}ds \right] e^{\int_{t_0}^{t}\rho(s)ds}
 \label{solution1}
\end{equation}
for $t \in I$, where $t_0 \in I$, $\bm{q}_0=\bm{q}(t_0)=\bm{\xi}(t_0)-\bm{p}(t_0)$ and $\bm{q}'_0=\bm{q}'(t_0)=\bm{\xi}'(t_0)-\bm{p}'(t_0)$. Hereafter, the proofs are given for each of the three cases (i)--(iii). 

Case (i). Assume that $f_{1}(t)$ and $f_{2}(t)$ given by \eqref{function1} and \eqref{function2} exist on $I$, and $\sup_{t \in I}f_{1}(t) < \infty$ and $\sup_{t \in I}f_{2}(t) < \infty$ are satisfied. Now we define
\begin{equation}
 \bm{c}_{1}:= \bm{q}'_0-\rho(t_0)\bm{q}_0+\int_{t_0}^{\sigma}\frac{e^{\int_{t_0}^{\mu}\left(\rho(\nu)+\frac{\beta(\nu)}{\alpha(\nu)}\right)d\nu}}{\alpha(\mu)}\bm{g}(\mu)d\mu.
 \label{specific_const1}
\end{equation}
Note that the integral contained in the right-hand side always converge. Actually, by using  $\sup_{t \in I}f_{1}(t) < \infty$, we can check that
\[ \left\|\int_{t_0}^{\sigma}\frac{e^{\int_{t_0}^{\mu}\left(\rho(\nu)+\frac{\beta(\nu)}{\alpha(\nu)}\right)d\nu}}{\alpha(\mu)}\bm{g}(\mu)d\mu\right\|
  \le \int_{t_0}^{\sigma}\frac{e^{\int_{t_0}^{\mu}\Re\left(\rho(\nu)+\frac{\beta(\nu)}{\alpha(\nu)}\right)d\nu}}{|\alpha(\mu)|}\|\bm{g}(\mu)\|d\mu
  \le \varepsilon f_{1}(t_0) < \infty. \]
That is, $\bm{c}_{1}$ is a well-defined constant vector. Therefore, \eqref{solution1} can be rewritten as
\begin{align}
 \bm{q}(t) &= \left[ \bm{q}_0+\int_{t_0}^{t} \left(\bm{c}_{1}-\int_{s}^{\sigma}\frac{e^{\int_{t_0}^{\mu}\left(\rho(\nu)+\frac{\beta(\nu)}{\alpha(\nu)}\right)d\nu}}{\alpha(\mu)}\bm{g}(\mu)d\mu \right) e^{-\int_{t_0}^{s}\left(2\rho(\mu)+\frac{\beta(\mu)}{\alpha(\mu)}\right)d\mu}ds \right] e^{\int_{t_0}^{t}\rho(s)ds} \nonumber\\
  &= \Bigg[ \bm{q}_0+\bm{c}_{1}\int_{t_0}^{t} e^{-\int_{t_0}^{s}\left(2\rho(\mu)+\frac{\beta(\mu)}{\alpha(\mu)}\right)d\mu}ds \nonumber\\
  & \hspace{5mm}-\int_{t_0}^{t}\left(\int_{s}^{\sigma}\frac{e^{\int_{s}^{\mu}\left(\rho(\nu)+\frac{\beta(\nu)}{\alpha(\nu)}\right)d\nu}}{\alpha(\mu)}\bm{g}(\mu)d\mu \right)e^{-\int_{t_0}^{s}\rho(\mu)d\mu}ds \Bigg] e^{\int_{t_0}^{t}\rho(s)ds}
 \label{solution2}
\end{align}
for $t \in I$. Moreover, we define
\[ \bm{c}_{2}:= \bm{q}_0-\int_{t_0}^{\sigma} \left(\int_{s}^{\sigma}\frac{e^{\int_{s}^{\mu}\left(\rho(\nu)+\frac{\beta(\nu)}{\alpha(\nu)}\right)d\nu}}{\alpha(\mu)}\bm{g}(\mu)d\mu \right)e^{-\int_{t_0}^{s}\rho(\mu)d\mu}ds. \]
Since
\begin{align*}
 &\left\|\int_{t_0}^{\sigma} \left(\int_{s}^{\sigma}\frac{e^{\int_{s}^{\mu}\left(\rho(\nu)+\frac{\beta(\nu)}{\alpha(\nu)}\right)d\nu}}{\alpha(\mu)}\bm{g}(\mu)d\mu \right)e^{-\int_{t_0}^{s}\rho(\mu)d\mu}ds\right\| \\
  &\le \int_{t_0}^{\sigma} \left\|\int_{s}^{\sigma}\frac{e^{\int_{s}^{\mu}\left(\rho(\nu)+\frac{\beta(\nu)}{\alpha(\nu)}\right)d\nu}}{\alpha(\mu)}\bm{g}(\mu)d\mu \right\| e^{-\int_{t_0}^{s}\Re(\rho(\mu))d\mu}ds\\
  &\le \int_{t_0}^{\sigma} \left(\int_{s}^{\sigma}\frac{e^{\int_{s}^{\mu}\Re\left(\rho(\nu)+\frac{\beta(\nu)}{\alpha(\nu)}\right)d\nu}}{|\alpha(\mu)|}\|\bm{g}(\mu)\|d\mu \right) e^{-\int_{t_0}^{s}\Re(\rho(\mu))d\mu}ds\\
  &\le \varepsilon\int_{t_0}^{\sigma} \left(\int_{s}^{\sigma}\frac{e^{\int_{s}^{\mu}\Re\left(\rho(\nu)+\frac{\beta(\nu)}{\alpha(\nu)}\right)d\nu}}{|\alpha(\mu)|}d\mu \right) e^{-\int_{t_0}^{s}\Re(\rho(\mu))d\mu}ds \\
  &\le \varepsilon \left(\sup_{t \in I}f_{1}(t)\right) \int_{t_0}^{\sigma} e^{-\int_{t_0}^{s}\Re(\rho(\mu))d\mu}ds \le \varepsilon \left(\sup_{t \in I}f_{1}(t)\right) f_{2}(t_0) < \infty
\end{align*}
holds, we can conclude that $\bm{c}_{2}$ is well-defined. Hence, \eqref{solution2} is rewritten as
\begin{align*}
 \bm{q}(t) &= \left[ \bm{c}_{2}+\bm{c}_{1}\int_{t_0}^{t} e^{-\int_{t_0}^{s}\left(2\rho(\mu)+\frac{\beta(\mu)}{\alpha(\mu)}\right)d\mu}ds+\int_{t}^{\sigma}\left(\int_{s}^{\sigma}\frac{e^{\int_{s}^{\mu}\left(\rho(\nu)+\frac{\beta(\nu)}{\alpha(\nu)}\right)d\nu}}{\alpha(\mu)}\bm{g}(\mu)d\mu \right)e^{-\int_{t_0}^{s}\rho(\mu)d\mu}ds \right] e^{\int_{t_0}^{t}\rho(s)ds}\\
  &= \left( \bm{c}_{2}+\bm{c}_{1}\int_{t_0}^{t} e^{-\int_{t_0}^{s}\left(2\rho(\mu)+\frac{\beta(\mu)}{\alpha(\mu)}\right)d\mu}ds\right) e^{\int_{t_0}^{t}\rho(s)ds} +\int_{t}^{\sigma}\left(\int_{s}^{\sigma}\frac{e^{\int_{s}^{\mu}\left(\rho(\nu)+\frac{\beta(\nu)}{\alpha(\nu)}\right)d\nu}}{\alpha(\mu)}\bm{g}(\mu)d\mu \right)e^{-\int_{t}^{s}\rho(\mu)d\mu}ds
\end{align*}
for $t \in I$. 

Next we consider the function
\[ \bm{w}(t) := \left( \bm{c}_{2}+\bm{c}_{1}\int_{t_0}^{t} e^{-\int_{t_0}^{s}\left(2\rho(\mu)+\frac{\beta(\mu)}{\alpha(\mu)}\right)d\mu}ds\right) e^{\int_{t_0}^{t}\rho(s)ds} \]
for $t \in I$. Note that this function is for $\bm{g}(t)\equiv \bm{0}$ in $\bm{q}(t)$ above. So it is a solution of the differential equation $\alpha(t)\bm{w}'' + \beta(t)\bm{w}' + \gamma(t)\bm{w} = \bm{0}$. Hence we see that the function
\[ \bm{x}_{1}(t) := \bm{w}(t)+\bm{p}(t) \]
is a solution of \eqref{oscillator1} for $t \in I$, where $\bm{p}(t)$ is a solution of \eqref{oscillator1} given at the beginning of the proof. This says that
\[ \bm{\xi}(t)-\bm{x}_{1}(t) = \bm{q}(t)-\bm{w}(t) = \int_{t}^{\sigma}\left(\int_{s}^{\sigma}\frac{e^{\int_{s}^{\mu}\left(\rho(\nu)+\frac{\beta(\nu)}{\alpha(\nu)}\right)d\nu}}{\alpha(\mu)}\bm{g}(\mu)d\mu \right)e^{-\int_{t}^{s}\rho(\mu)d\mu}ds, \]
and so that
\[ \|\bm{\xi}(t)-\bm{x}_{1}(t)\| \le \varepsilon \int_{t}^{\sigma}\left(\int_{s}^{\sigma}\frac{e^{\int_{s}^{\mu}\Re\left(\rho(\nu)+\frac{\beta(\nu)}{\alpha(\nu)}\right)d\nu}}{|\alpha(\mu)|}d\mu \right)e^{-\int_{t}^{s}\Re(\rho(\mu))d\mu}ds 
  \le \varepsilon \left(\sup_{t \in I}f_{1}(t)\right)\left(\sup_{t \in I}f_{2}(t)\right) < \infty \]
for $t \in I$. Thus, \eqref{oscillator1} is Ulam stable on $I$. Moreover, 
\[ L_{1} = \sup_{t \in I} \int_{t}^{\sigma} \left(\int_{s}^{\sigma} \frac{e^{\int_s^\mu \Re\left(\rho(\nu)+\frac{\beta(\nu)}{\alpha(\nu)}\right) d\nu}}{|\alpha(\mu)|}d\mu\right) e^{-\int_t^s \Re(\rho(\mu)) d\mu}ds \]
is an Ulam constant for \eqref{oscillator1}. 

Case (ii). Assume that $f_{1}(t)$ and $f_{3}(t)$ given by \eqref{function1} and \eqref{function3} exist on $I$, and $\sup_{t \in I}f_{1}(t) < \infty$ and $\sup_{t \in I}f_{3}(t) < \infty$ hold. As in Case (i), we can rewrite $\bm{q}(t)$ as \eqref{solution2}, using the constant $\bm{c}_{1}$ defined in \eqref{specific_const1}. Now we define
\[ \bm{c}_{3}:= \bm{q}_0+\int_{\tau}^{t_0} \left(\int_{s}^{\sigma}\frac{e^{\int_{s}^{\mu}\left(\rho(\nu)+\frac{\beta(\nu)}{\alpha(\nu)}\right)d\nu}}{\alpha(\mu)}\bm{g}(\mu)d\mu \right)e^{\int_{s}^{t_0}\rho(\mu)d\mu}ds. \]
Since
\begin{align*}
  &\left\|\int_{\tau}^{t_0} \left(\int_{s}^{\sigma}\frac{e^{\int_{s}^{\mu}\left(\rho(\nu)+\frac{\beta(\nu)}{\alpha(\nu)}\right)d\nu}}{\alpha(\mu)}\bm{g}(\mu)d\mu \right)e^{\int_{s}^{t_0}\rho(\mu)d\mu}ds\right\| \\
  &\le \int_{\tau}^{t_0} \left(\int_{s}^{\sigma}\frac{e^{\int_{s}^{\mu}\Re\left(\rho(\nu)+\frac{\beta(\nu)}{\alpha(\nu)}\right)d\nu}}{|\alpha(\mu)|}\|\bm{g}(\mu)\|d\mu \right) e^{\int_{s}^{t_0}\Re(\rho(\mu))d\mu}ds\\
  &\le \varepsilon \left(\sup_{t \in I}f_{1}(t)\right) \int_{\tau}^{t_0} e^{\int_{s}^{t_0}\Re(\rho(\mu))d\mu}ds \le \varepsilon \left(\sup_{t \in I}f_{1}(t)\right) f_{3}(t_0) < \infty
\end{align*}
holds, $\bm{c}_{3}$ is a well-defined constant. Thus, \eqref{solution2} is rewritten as
\begin{align*}
 \bm{q}(t) &= \left[ \bm{c}_{3}+\bm{c}_{1}\int_{t_0}^{t} e^{-\int_{t_0}^{s}\left(2\rho(\mu)+\frac{\beta(\mu)}{\alpha(\mu)}\right)d\mu}ds-\int_{\tau}^{t}\left(\int_{s}^{\sigma}\frac{e^{\int_{s}^{\mu}\left(\rho(\nu)+\frac{\beta(\nu)}{\alpha(\nu)}\right)d\nu}}{\alpha(\mu)}\bm{g}(\mu)d\mu \right)e^{\int_{s}^{t_0}\rho(\mu)d\mu}ds \right] e^{\int_{t_0}^{t}\rho(s)ds}\\
  &= \left( \bm{c}_{3}+\bm{c}_{1}\int_{t_0}^{t} e^{-\int_{t_0}^{s}\left(2\rho(\mu)+\frac{\beta(\mu)}{\alpha(\mu)}\right)d\mu}ds\right) e^{\int_{t_0}^{t}\rho(s)ds} -\int_{\tau}^{t}\left(\int_{s}^{\sigma}\frac{e^{\int_{s}^{\mu}\left(\rho(\nu)+\frac{\beta(\nu)}{\alpha(\nu)}\right)d\nu}}{\alpha(\mu)}\bm{g}(\mu)d\mu \right)e^{\int_{s}^{t}\rho(\mu)d\mu}ds
\end{align*}
for $t \in I$. 

Next we consider the solution of \eqref{oscillator1} given by
\[ \bm{x}_{2}(t) := \left( \bm{c}_{3}+\bm{c}_{1}\int_{t_0}^{t} e^{-\int_{t_0}^{s}\left(2\rho(\mu)+\frac{\beta(\mu)}{\alpha(\mu)}\right)d\mu}ds\right) e^{\int_{t_0}^{t}\rho(s)ds}+\bm{p}(t) \]
for $t \in I$. Then
\[ \bm{x}_{2}(t)-\bm{\xi}(t) = \int_{\tau}^{t}\left(\int_{s}^{\sigma}\frac{e^{\int_{s}^{\mu}\left(\rho(\nu)+\frac{\beta(\nu)}{\alpha(\nu)}\right)d\nu}}{\alpha(\mu)}\bm{g}(\mu)d\mu \right)e^{\int_{s}^{t}\rho(\mu)d\mu}ds, \]
and so that
\[ \|\bm{x}_{2}(t)-\bm{\xi}(t)\| \le \varepsilon \int_{\tau}^{t}\left(\int_{s}^{\sigma}\frac{e^{\int_{s}^{\mu}\Re\left(\rho(\nu)+\frac{\beta(\nu)}{\alpha(\nu)}\right)d\nu}}{|\alpha(\mu)|}d\mu \right)e^{\int_{s}^{t}\Re(\rho(\mu))d\mu}ds 
  \le \varepsilon \left(\sup_{t \in I}f_{1}(t)\right)\left(\sup_{t \in I}f_{3}(t)\right) < \infty \]
for $t \in I$. Thus, \eqref{oscillator1} is Ulam stable on $I$. Moreover, 
\[ L_{2} = \sup_{t \in I} \int_{\tau}^{t} \left(\int_{s}^{\sigma} \frac{e^{\int_s^\mu \Re\left(\rho(\nu)+\frac{\beta(\nu)}{\alpha(\nu)}\right) d\nu}}{|\alpha(\mu)|}d\mu\right) e^{\int_s^t \Re(\rho(\mu)) d\mu}ds \]
is an Ulam constant for \eqref{oscillator1}. 

Case (iii). Assume that $f_{3}(t)$ and $f_{4}(t)$ given by \eqref{function3} and \eqref{function4} exist on $I$, and $\sup_{t \in I}f_{3}(t) < \infty$ and $\sup_{t \in I}f_{4}(t) < \infty$ are satisfied. Now we define
\begin{equation*}
 \bm{c}_{4}:= \bm{q}'_0-\rho(t_0)\bm{q}_0-\int_{\tau}^{t_0}\frac{e^{-\int_{\mu}^{t_0}\left(\rho(\nu)+\frac{\beta(\nu)}{\alpha(\nu)}\right)d\nu}}{\alpha(\mu)}\bm{g}(\mu)d\mu.
 \label{specific_const4}
\end{equation*}
Since  $\sup_{t \in I}f_{4}(t) < \infty$ holds, we see that
\[ \left\|\int_{\tau}^{t_0}\frac{e^{-\int_{\mu}^{t_0}\left(\rho(\nu)+\frac{\beta(\nu)}{\alpha(\nu)}\right)d\nu}}{\alpha(\mu)}\bm{g}(\mu)d\mu\right\|
  \le \int_{\tau}^{t_0}\frac{e^{-\int_{\mu}^{t_0}\Re\left(\rho(\nu)+\frac{\beta(\nu)}{\alpha(\nu)}\right)d\nu}}{|\alpha(\mu)|}\|\bm{g}(\mu)\|d\mu
  \le \varepsilon f_{4}(t_0) < \infty, \]
and that $\bm{c}_{4}$ is a well-defined constant vector. Therefore, \eqref{solution1} can be rewritten as
\[ \bm{q}(t) = \Bigg[ \bm{q}_0+\bm{c}_{4}\int_{t_0}^{t} e^{-\int_{t_0}^{s}\left(2\rho(\mu)+\frac{\beta(\mu)}{\alpha(\mu)}\right)d\mu}ds + \int_{t_0}^{t}\left(\int_{\tau}^{s}\frac{e^{-\int_{\mu}^{s}\left(\rho(\nu)+\frac{\beta(\nu)}{\alpha(\nu)}\right)d\nu}}{\alpha(\mu)}\bm{g}(\mu)d\mu \right)e^{-\int_{t_0}^{s}\rho(\mu)d\mu}ds \Bigg] e^{\int_{t_0}^{t}\rho(s)ds} \]
for $t \in I$. Moreover, we define
\[ \bm{c}_{5}:= \bm{q}_0-\int_{\tau}^{t_0} \left(\int_{\tau}^{s}\frac{e^{-\int_{\mu}^{s}\left(\rho(\nu)+\frac{\beta(\nu)}{\alpha(\nu)}\right)d\nu}}{\alpha(\mu)}\bm{g}(\mu)d\mu \right)e^{\int_{s}^{t_0}\rho(\mu)d\mu}ds. \]
Since
\begin{align*}
 &\left\|\int_{\tau}^{t_0} \left(\int_{\tau}^{s}\frac{e^{-\int_{\mu}^{s}\left(\rho(\nu)+\frac{\beta(\nu)}{\alpha(\nu)}\right)d\nu}}{\alpha(\mu)}\bm{g}(\mu)d\mu \right)e^{\int_{s}^{t_0}\rho(\mu)d\mu}ds\right\| \\
  &\le \int_{\tau}^{t_0} \left(\int_{\tau}^{s}\frac{e^{-\int_{\mu}^{s}\Re\left(\rho(\nu)+\frac{\beta(\nu)}{\alpha(\nu)}\right)d\nu}}{|\alpha(\mu)|}\|\bm{g}(\mu)\|d\mu \right) e^{\int_{s}^{t_0}\Re(\rho(\mu))d\mu}ds\\
  &\le \varepsilon \left(\sup_{t \in I}f_{4}(t)\right) \int_{\tau}^{t_0} e^{\int_{s}^{t_0}\Re(\rho(\mu))d\mu}ds \le \varepsilon \left(\sup_{t \in I}f_{4}(t)\right) f_{3}(t_0) < \infty
\end{align*}
holds, we see that $\bm{c}_{5}$ is well-defined. Hence, $\bm{q}(t)$ is rewritten as
\[ \bm{q}(t) = \left( \bm{c}_{5}+\bm{c}_{4}\int_{t_0}^{t} e^{-\int_{t_0}^{s}\left(2\rho(\mu)+\frac{\beta(\mu)}{\alpha(\mu)}\right)d\mu}ds\right) e^{\int_{t_0}^{t}\rho(s)ds} +\int_{\tau}^{t}\left(\int_{\tau}^{s}\frac{e^{-\int_{\mu}^{s}\left(\rho(\nu)+\frac{\beta(\nu)}{\alpha(\nu)}\right)d\nu}}{\alpha(\mu)}\bm{g}(\mu)d\mu \right)e^{\int_{s}^{t}\rho(\mu)d\mu}ds \]
for $t \in I$. 

Next we consider the solution of \eqref{oscillator1} given by
\[ \bm{x}_{3}(t) := \left( \bm{c}_{5}+\bm{c}_{4}\int_{t_0}^{t} e^{-\int_{t_0}^{s}\left(2\rho(\mu)+\frac{\beta(\mu)}{\alpha(\mu)}\right)d\mu}ds\right) e^{\int_{t_0}^{t}\rho(s)ds}+\bm{p}(t) \]
for $t \in I$. Then
\[ \bm{\xi}(t)-\bm{x}_{3}(t) = \int_{\tau}^{t}\left(\int_{\tau}^{s}\frac{e^{-\int_{\mu}^{s}\left(\rho(\nu)+\frac{\beta(\nu)}{\alpha(\nu)}\right)d\nu}}{\alpha(\mu)}\bm{g}(\mu)d\mu \right)e^{\int_{s}^{t}\rho(\mu)d\mu}ds, \]
and so that
\[ \|\bm{\xi}(t)-\bm{x}_{3}(t)\| \le \varepsilon \int_{\tau}^{t}\left(\int_{\tau}^{s}\frac{e^{-\int_{\mu}^{s}\Re\left(\rho(\nu)+\frac{\beta(\nu)}{\alpha(\nu)}\right)d\nu}}{|\alpha(\mu)|}d\mu \right)e^{\int_{s}^{t}\Re(\rho(\mu))d\mu}ds 
  \le \varepsilon \left(\sup_{t \in I}f_{3}(t)\right)\left(\sup_{t \in I}f_{4}(t)\right) < \infty \]
for $t \in I$. Thus, \eqref{oscillator1} is Ulam stable on $I$. Moreover, 
\[ L_{3} = \sup_{t \in I} \int_{\tau}^{t}\left(\int_{\tau}^{s}\frac{e^{-\int_{\mu}^{s}\Re\left(\rho(\nu)+\frac{\beta(\nu)}{\alpha(\nu)}\right)d\nu}}{|\alpha(\mu)|}d\mu \right)e^{\int_{s}^{t}\Re(\rho(\mu))d\mu}ds \]
is an Ulam constant for \eqref{oscillator1}. The proof is now complete.
\end{proof}

Consider the constant coefficients linear oscillator
\begin{equation}
 a_0 \bm{x}'' + a_1 \bm{x}' + a_2 \bm{x} = 0,
 \label{constoscillator}
\end{equation}
where $a_0$, $a_1$ and $a_2$ are complex-valued constants, and $a_0 \ne 0$. Then we obtain the following result. 


\begin{corollary}\label{constcoro1}
Let $I=\R$, and $\lambda_1$ and $\lambda_2$ be the roots of the characteristic equation
\begin{equation*}
 a_0 \lambda^2 + a_1 \lambda + a_2 = 0.
 \label{characteristic}
\end{equation*}
If $a_0 \Re(\lambda_1)\Re(\lambda_2) \ne 0$, then \eqref{constoscillator} is Ulam stable on $\R$, and an Ulam constant is $\frac{1}{|a_0\Re(\lambda_1)\Re(\lambda_2)|}$. 
\end{corollary}

\begin{proof}
Assume $a_0\ne0$. The proof is divided into three cases (i) $\Re(\lambda_1) \ge \Re(\lambda_2) > 0$, (ii) $\Re(\lambda_1) > 0 > \Re(\lambda_2)$, and (iii) $0 > \Re(\lambda_1) \ge \Re(\lambda_2)$. First, we notice that, since $\lambda_1$ and $\lambda_2$ are the roots of the characteristic equation, they are constant solutions of \eqref{Riccati}; that is, we can choose $\rho(t)= \lambda_1$ or $\rho(t)= \lambda_2$ for all $t \in \R$. Moreover, if
\[ \rho(t) = \lambda_2 = \frac{-a_1+\sqrt{a_1^2-4a_0a_2}}{2 a_0}, \]
then
\[ \rho(t)+\frac{\beta(t)}{\alpha(t)} = \lambda_2+\frac{a_1}{a_0} = -\lambda_1. \]
We will use Theorem \ref{UlamS} with $\tau=-\infty$ and $\sigma=\infty$. 

Case (i). Suppose $\Re(\lambda_1) \ge \Re(\lambda_2) > 0$. In this case, Theorem \ref{UlamS} (i) will be used. Set $\rho(t)= \lambda_2$. From
\[f_{1}(t) = \int_{t}^{\sigma} \frac{e^{\int_t^s \Re\left(\rho(\mu)+\frac{\beta(\mu)}{\alpha(\mu)}\right) d\mu}}{|\alpha(s)|}ds = \int_{t}^{\infty} \frac{e^{-\int_t^s \Re(\lambda_1) d\mu}}{|a_0|}ds = \frac{1}{|a_0|\Re(\lambda_1)} \]
and
\[ f_{2}(t) = \int_{t}^{\sigma} e^{-\int_t^s \Re(\rho(\mu)) d\mu}ds = \int_{t}^{\infty} e^{-\int_t^s \Re(\lambda_2) d\mu}ds = \frac{1}{\Re(\lambda_2)} \]
for all $t \in \R$, $f_{1}(t)$ and $f_{2}(t)$ exist and are bounded on $\R$. Hence, by Theorem \ref{UlamS} (i), \eqref{constoscillator} is Ulam stable on $\R$, and an Ulam constant is
\[ \sup_{t \in I} \int_{t}^{\sigma} \left(\int_{s}^{\sigma} \frac{e^{\int_s^\mu \Re\left(\rho(\nu)+\frac{\beta(\nu)}{\alpha(\nu)}\right) d\nu}}{|\alpha(\mu)|}d\mu\right) e^{-\int_t^s \Re(\rho(\mu)) d\mu}ds = \sup_{t \in I} f_{1}(t)f_{2}(t) = \frac{1}{|a_0|\Re(\lambda_1) \Re(\lambda_2)}. \]

Case (ii). Suppose $\Re(\lambda_1) > 0 > \Re(\lambda_2)$. Set $\rho(t)= \lambda_2$. From $f_{1}(t) = \frac{1}{|a_0|\Re(\lambda_1)}$ and
\[ f_{3}(t) = \int_{\tau}^{t} e^{\int_s^t \Re(\rho(\mu)) d\mu}ds = \int_{-\infty}^{t} e^{\int_s^t \Re(\lambda_2) d\mu}ds = \frac{1}{-\Re(\lambda_2)} \]
for all $t \in \R$, $f_{1}(t)$ and $f_{3}(t)$ exist and are bounded on $\R$. Hence, by Theorem \ref{UlamS} (ii), \eqref{constoscillator} is Ulam stable on $\R$, and an Ulam constant is
\[ \sup_{t \in I} \int_{\tau}^{t} \left(\int_{s}^{\sigma} \frac{e^{\int_s^\mu \Re\left(\rho(\nu)+\frac{\beta(\nu)}{\alpha(\nu)}\right) d\nu}}{|\alpha(\mu)|}d\mu\right) e^{\int_s^t \Re(\rho(\mu)) d\mu}ds = \sup_{t \in I} f_{1}(t)f_{3}(t) = \frac{1}{|a_0\Re(\lambda_1) \Re(\lambda_2)|}. \]

Case (iii). Suppose $0 > \Re(\lambda_1) \ge \Re(\lambda_2)$. Set $\rho(t)= \lambda_2$. From $f_{3}(t) = \frac{1}{-\Re(\lambda_2)}$ and
\[ f_{4}(t) = \int_{\tau}^{t} \frac{e^{-\int_s^t \Re\left(\rho(\mu)+\frac{\beta(\mu)}{\alpha(\mu)}\right) d\mu}}{|\alpha(s)|}ds = \int_{-\infty}^{t} \frac{e^{\int_s^t \Re(\lambda_1) d\mu}}{|a_0|}ds = \frac{1}{-|a_0|\Re(\lambda_1)} \]
for all $t \in \R$, $f_{3}(t)$ and $f_{4}(t)$ exist and are bounded on $\R$. Hence, by Theorem \ref{UlamS} (iii), \eqref{constoscillator} is Ulam stable on $\R$, and an Ulam constant is
\[ \sup_{t \in I} \int_{\tau}^{t}\left(\int_{\tau}^{s}\frac{e^{-\int_{\mu}^{s}\Re\left(\rho(\nu)+\frac{\beta(\nu)}{\alpha(\nu)}\right)d\nu}}{|\alpha(\mu)|}d\mu \right)e^{\int_{s}^{t}\Re(\rho(\mu))d\mu}ds = \sup_{t \in I} f_{3}(t)f_{4}(t) = \frac{1}{|a_0\Re(\lambda_1) \Re(\lambda_2)|}. \]
Therefore, in any case, an Ulam constant is $\frac{1}{|a_0\Re(\lambda_1) \Re(\lambda_2)|}$.
\end{proof}


\section{Minimum Ulam constants}

In this section, we show that the Ulam constants given in Theorem \ref{UlamS} are the minimum Ulam constants by restricting to real-valued scalar functions.


\begin{theorem}\label{thm:minimum1}
Let $I$ be either $(\tau,\sigma)$, $(\tau,\sigma]$, $[\tau,\sigma)$ or $[\tau,\sigma]$, where $-\infty \le \tau < \sigma \le \infty$. Suppose that $\alpha$, $\beta$, $\gamma: I \to \R$ are real-valued continuous functions, and $\alpha(t)\ne 0$ for all $t \in I$, and there exists a real-valued solution $\rho: I \to \R$ of \eqref{Riccati}. Then the following (i), (ii) and (iii) below hold:
\begin{itemize}
  \item[(i)] suppose that $f_{1}(t)$ and $f_{2}(t)$ given by \eqref{function1} and \eqref{function2} exist for all $t\in I$, and $\sup_{t \in I}f_{1}(t) < \infty$ and $\sup_{t \in I}f_{2}(t) < \infty$ hold. If
\begin{equation}
 \lim_{t\to \sigma^{-}} \int_{t_0}^t \rho(s) ds = \infty \quad\text{and}\quad \lim_{t\to \sigma^{-}} e^{\int_{t_0}^{t}\rho(s)ds} \int_{t_0}^{t} e^{-\int_{t_0}^{s}\left(2\rho(\mu)+\frac{\beta(\mu)}{\alpha(\mu)}\right)d\mu}ds = \infty, \quad t_0\in(\tau,\sigma),
 \label{diverge1}
\end{equation}
then \eqref{oscillator1} is Ulam stable on $I$, and the minimum Ulam constant is
\begin{equation}
 B_{1} := \sup_{t \in I} \int_{t}^{\sigma} \left(\int_{s}^{\sigma} \frac{e^{\int_s^\mu \left(\rho(\nu)+\frac{\beta(\nu)}{\alpha(\nu)}\right) d\nu}}{|\alpha(\mu)|}d\mu\right) e^{-\int_t^s \rho(\mu) d\mu}ds;
 \label{minimum1}
\end{equation}
  \item[(ii)] suppose that $f_{1}(t)$ and $f_{3}(t)$ given by \eqref{function1} and \eqref{function3} exist for all $t\in I$, and $\sup_{t \in I}f_{1}(t) < \infty$ and $\sup_{t \in I}f_{3}(t) < \infty$ hold. If
\begin{equation}
 \lim_{t\to \tau^{+}} \int_{t_0}^t \rho(s) ds = \infty \quad\text{and}\quad \lim_{t\to \sigma^{-}} e^{\int_{t_0}^{t}\rho(s)ds} \int_{t_0}^{t} e^{-\int_{t_0}^{s}\left(2\rho(\mu)+\frac{\beta(\mu)}{\alpha(\mu)}\right)d\mu}ds = \infty, \quad t_0\in(\tau,\sigma),
 \label{diverge2}
\end{equation}
then \eqref{oscillator1} is Ulam stable on $I$, and the minimum Ulam constant is
\begin{equation*}
 B_{2} := \sup_{t \in I} \int_{\tau}^{t} \left(\int_{s}^{\sigma} \frac{e^{\int_s^\mu \left(\rho(\nu)+\frac{\beta(\nu)}{\alpha(\nu)}\right) d\nu}}{|\alpha(\mu)|}d\mu\right) e^{\int_s^t \rho(\mu) d\mu}ds;
 \label{minimum2}
\end{equation*}
  \item[(iii)] suppose that $f_{3}(t)$ and $f_{4}(t)$ given by \eqref{function3} and \eqref{function4} exist for all $t\in I$, and $\sup_{t \in I}f_{3}(t) < \infty$ and $\sup_{t \in I}f_{4}(t) < \infty$ hold. If
\begin{equation}
 \lim_{t\to \tau^{+}} \int_{t_0}^t \rho(s) ds = \infty \quad\text{and}\quad \lim_{t\to \tau^{+}} e^{\int_{t_0}^{t}\rho(s)ds} \int_{t_0}^{t} e^{-\int_{t_0}^{s}\left(2\rho(\mu)+\frac{\beta(\mu)}{\alpha(\mu)}\right)d\mu}ds = -\infty, \quad t_0\in(\tau,\sigma),
 \label{diverge3}
\end{equation}
then \eqref{oscillator1} is Ulam stable on $I$, and the minimum Ulam constant is
\begin{equation}
 B_{3} := \sup_{t \in I} \int_{\tau}^{t} \left(\int_{\tau}^{s} \frac{e^{-\int_\mu^s \left(\rho(\nu)+\frac{\beta(\nu)}{\alpha(\nu)}\right) d\nu}}{|\alpha(\mu)|}d\mu\right) e^{\int_s^t \rho(\mu) d\mu}ds.
 \label{minimum3}
\end{equation}
\end{itemize}
\end{theorem}

\begin{proof}
Assume that $\alpha$, $\beta$ and $\gamma$ are real-valued continuous functions, and $\alpha(t)\ne 0$ for all $t \in I$. Let $\rho$ be a real-valued solution of \eqref{Riccati}. Throughout this proof, let $f_{1}$, $f_{2}$, $f_{3}$ and $f_{4}$ be the functions defined by \eqref{function1}--\eqref{function4}, respectively. Note that $f_{1}$, $f_{2}$, $f_{3}$ and $f_{4}$ are real-valued functions on $I$. Let $t_0\in(\tau,\sigma)$. 

Case (i). Assume that $f_{1}(t)$ and $f_{2}(t)$ exist for all $t\in I$, and $\sup_{t \in I}f_{1}(t) < \infty$ and $\sup_{t \in I}f_{2}(t) < \infty$ hold. By Theorem \ref{UlamS}, we see that \eqref{oscillator1} is Ulam stable on $I$, with an Ulam constant $B_1$, where $B_1$ is defined by \eqref{minimum1}. Let $\varepsilon>0$. Now we consider the function
\[ \bm{q}(t):= \left( \bm{c}_{2}+\bm{c}_{1}\int_{t_0}^{t} e^{-\int_{t_0}^{s}\left(2\rho(\mu)+\frac{\beta(\mu)}{\alpha(\mu)}\right)d\mu}ds\right) e^{\int_{t_0}^{t}\rho(s)ds} + \varepsilon \left[\int_{t}^{\sigma}\left(\int_{s}^{\sigma}\frac{e^{\int_{s}^{\mu}\left(\rho(\nu)+\frac{\beta(\nu)}{\alpha(\nu)}\right)d\nu}}{\alpha(\mu)}d\mu \right)e^{-\int_{t}^{s}\rho(\mu)d\mu}ds\right] \bm{u} \]
for all $t \in I$, where $\bm{c}_{1}$ and $\bm{c}_{2}$ are well-defined constants given in the proof of Theorem \ref{UlamS} (i) and $\bm{u}$ is the unit vector. From the proof of Theorem \ref{UlamS} (i), we find that $\bm{q}(t)$ is a solution of the equation
\[ \alpha(t)\bm{q}'' + \beta(t)\bm{q}' + \gamma(t)\bm{q} = \varepsilon \bm{u}. \]
Let $\bm{p}(t)$ be a solution to \eqref{oscillator1} on $I$, and let $\bm{\xi}(t):=\bm{q}(t)+\bm{p}(t)$ for $t \in I$. Then
\begin{equation}
 \|\alpha(t)\bm{\xi}'' + \beta(t)\bm{\xi}' + \gamma(t)\bm{\xi} - \bm{f}(t)\| = \varepsilon
 \label{eqallity1}
\end{equation}
is satisfied for $t \in I$. By the Ulam stability for \eqref{oscillator1}, we find that there exists a solution $\bm{x}_{1}: I \to \R^n$ of \eqref{oscillator1} such that
\begin{equation}
 \sup_{t \in I}\|\bm{\xi}(t)-\bm{x}_{1}(t)\| \le B_1 \varepsilon.
 \label{ineqallity1}
\end{equation}
More precisely, from the proof of Theorem \ref{UlamS} (i), we know that $\bm{x}_{1}(t)$ is given as
\[ \bm{x}_{1}(t) = \left( \bm{c}_{2}+\bm{c}_{1}\int_{t_0}^{t} e^{-\int_{t_0}^{s}\left(2\rho(\mu)+\frac{\beta(\mu)}{\alpha(\mu)}\right)d\mu}ds\right) e^{\int_{t_0}^{t}\rho(s)ds} + \bm{p}(t). \]
We show that $B_1$ is the minimum Ulam constant by using the following two steps.

Step 1. We first show that $\bm{x}_{1}(t)$ is the unique solution of \eqref{oscillator1} satisfying \eqref{ineqallity1}. To show this fact using contradiction, we assume that there exists a solution $\bm{y}_{1}(t)$ of \eqref{oscillator1} such that $\bm{y}_{1}(t) \neq \bm{x}_{1}(t)$ for all $t \in I$. That is, $\bm{y}_{1}(t)$ is written as
\[ \bm{y}_{1}(t) = \left( \bm{d}_{2}+\bm{d}_{1}\int_{t_0}^{t} e^{-\int_{t_0}^{s}\left(2\rho(\mu)+\frac{\beta(\mu)}{\alpha(\mu)}\right)d\mu}ds\right) e^{\int_{t_0}^{t}\rho(s)ds} + \bm{p}(t) \]
with $(\bm{d}_{1}, \bm{d}_{2}) \ne (\bm{c}_{1}, \bm{c}_{2})$. Thus, we have
\[ \|\bm{y}_{1}(t)-\bm{x}_{1}(t)\| \le \|\bm{\xi}(t)-\bm{y}_{1}(t)\| + \|\bm{\xi}(t)-\bm{x}_{1}(t)\| \le 2 B_1 \varepsilon \]
for all $t \in I$. However, with \eqref{diverge1}, the following holds: 
\[ \lim_{t\to\sigma^-}\|\bm{y}_{1}(t)-\bm{x}_{1}(t)\| 
  = \lim_{t\to\sigma^-}\left\|\left[ (\bm{d}_{2}-\bm{c}_{2})+(\bm{d}_{1}-\bm{c}_{1})\int_{t_0}^{t} e^{-\int_{t_0}^{s}\left(2\rho(\mu)+\frac{\beta(\mu)}{\alpha(\mu)}\right)d\mu}ds\right] e^{\int_{t_0}^{t}\rho(s)ds}\right\| = \infty. \]
This contradicts the above inequality. 

Step 2. We next show that $B_1$ is the minimum Ulam constant. By way of contradiction, we assume that there exists $0<U_1<B_1$ such that
\[ \sup_{t \in I}\|\bm{\xi}(t)-\bm{x}_{1}(t)\| \le U_1 \varepsilon. \]
Note that $\bm{\xi}(t)$ satisfies \eqref{eqallity1}, and that there is no other possible solution of \eqref{oscillator1} that satisfies this inequality other than $\bm{x}_{1}(t)$ by Step 1. 
However, we see that
\begin{align*}
 \sup_{t \in I}\|\bm{\xi}(t)-\bm{x}_{1}(t)\| &= \sup_{t \in I}\left\|\varepsilon \left[\int_{t}^{\sigma}\left(\int_{s}^{\sigma}\frac{e^{\int_{s}^{\mu}\left(\rho(\nu)+\frac{\beta(\nu)}{\alpha(\nu)}\right)d\nu}}{\alpha(\mu)}d\mu \right)e^{-\int_{t}^{s}\rho(\mu)d\mu}ds\right] \bm{u} \right\| = B_1 \varepsilon,
\end{align*}
and thus,
\[ \sup_{t \in I}\|\bm{\xi}(t)-\bm{x}_{1}(t)\| \le U_1 \varepsilon < B_1 \varepsilon = \sup_{t \in I}\|\bm{\xi}(t)-\bm{x}_{1}(t)\|. \]
This is a contradiction. Hence we can conclude that $B_1$ is the minimum Ulam constant. 

Cases (ii) and (iii) can be shown by the same technique as in Case (i). The proof is now complete.
\end{proof}

Consider the linear oscillator \eqref{constoscillator} again, where $a_0$, $a_1$ and $a_2$ are real-valued constants, and $a_0 \ne 0$. Theorem \ref{thm:minimum1} implies the following result. 


\begin{corollary}\label{constminimum1}
Let $I=\R$, and $\lambda_1$ and $\lambda_2$ be non-zero real roots of the characteristic equation
\[ a_0 \lambda^2 + a_1 \lambda + a_2 = 0. \]
If $a_0 \ne 0$, then \eqref{constoscillator} is Ulam stable on $\R$, and the minimum Ulam constant is $\frac{1}{|a_0\lambda_1\lambda_2|}$. 
\end{corollary}

\begin{proof}
The proof is divided into three cases (i) $\lambda_1 \ge \lambda_2 > 0$, (ii) $\lambda_1 > 0 > \lambda_2$, and (iii) $0 > \lambda_1 \ge \lambda_2$. Recall the proof of Corollary \ref{constcoro1}. Since $\lambda_1$ and $\lambda_2$ are the roots of the characteristic equation, they are constant solutions of \eqref{Riccati}, and if $\rho(t) = \lambda_2$, then $\rho(t)+\frac{\beta(t)}{\alpha(t)} = -\lambda_1$. We will use Theorem \ref{thm:minimum1} with $\tau=-\infty$ and $\sigma=\infty$. Let $t_0 \in (-\infty,\infty)$. 

Case (i). Suppose $\lambda_1 \ge \lambda_2 > 0$. Set $\rho(t) = \lambda_2$. From the facts obtained in Case (i) in the proof of Corollary \ref{constcoro1}, we have $f_{1}(t) = \frac{1}{|a_0|\lambda_1}$ and $f_{2}(t) = \frac{1}{\lambda_2}$, and thus, $f_{1}(t)$ and $f_{2}(t)$ exist and are bounded on $\R$. Moreover, 
\[ \lim_{t\to \sigma^{-}} \int_{t_0}^t \rho(s) ds = \lim_{t\to \infty} \int_{t_0}^t \lambda_2 ds = \infty \]
and
\[ \lim_{t\to \sigma^{-}} e^{\int_{t_0}^{t}\rho(s)ds} \int_{t_0}^{t} e^{-\int_{t_0}^{s}\left(2\rho(\mu)+\frac{\beta(\mu)}{\alpha(\mu)}\right)d\mu}ds 
  = \lim_{t\to \infty} e^{\int_{t_0}^{t}\lambda_2ds} \int_{t_0}^{t} e^{\int_{t_0}^{s}(\lambda_1-\lambda_2)d\mu}ds = \infty \]
are satisfied; that is, \eqref{diverge1} holds. Hence, by Theorem \ref{thm:minimum1} (i), \eqref{constoscillator} is Ulam stable on $I$, and the minimum Ulam constant is $\frac{1}{|a_0|\lambda_1\lambda_2}$. 

Case (ii). Suppose $\lambda_1 > 0 > \lambda_2$. Set $\rho(t) = \lambda_2$. From the facts obtained in Case (ii) in the proof of Corollary \ref{constcoro1}, we have $f_{1}(t) = \frac{1}{|a_0|\lambda_1}$ and $f_{3}(t) = \frac{1}{-\lambda_2}$, and thus, $f_{1}(t)$ and $f_{3}(t)$ exist and are bounded on $\R$. Moreover, 
\[ \lim_{t\to \tau^{+}} \int_{t_0}^t \rho(s) ds = \lim_{t\to -\infty} \int_{t_0}^t \lambda_2 ds = \infty \]
and
\[ \lim_{t\to \sigma^{-}} e^{\int_{t_0}^{t}\rho(s)ds} \int_{t_0}^{t} e^{-\int_{t_0}^{s}\left(2\rho(\mu)+\frac{\beta(\mu)}{\alpha(\mu)}\right)d\mu}ds 
  = \lim_{t\to \infty} e^{\int_{t_0}^{t}\lambda_2ds} \int_{t_0}^{t} e^{\int_{t_0}^{s}(\lambda_1-\lambda_2)d\mu}ds = \infty \]
are satisfied; that is, \eqref{diverge2} holds. Hence, by Theorem \ref{thm:minimum1} (ii), \eqref{constoscillator} is Ulam stable on $I$, and the minimum Ulam constant is $\frac{1}{|a_0\lambda_1\lambda_2|}$. 

Case (iii). Suppose $0 > \lambda_1 \ge \lambda_2$. Set $\rho(t) = \lambda_2$. From the facts obtained in Case (iii) in the proof of Corollary \ref{constcoro1}, we have $f_{3}(t) = \frac{1}{-\lambda_2}$ and $f_{4}(t) = \frac{1}{-|a_0|\lambda_1}$, and thus, $f_{3}(t)$ and $f_{4}(t)$ exist and are bounded on $\R$. Moreover, $\lim_{t\to \tau^{+}} \int_{t_0}^t \rho(s) ds = \infty$ and
\[ \lim_{t\to \tau^{+}} e^{\int_{t_0}^{t}\rho(s)ds} \int_{t_0}^{t} e^{-\int_{t_0}^{s}\left(2\rho(\mu)+\frac{\beta(\mu)}{\alpha(\mu)}\right)d\mu}ds 
  = \lim_{t\to -\infty} e^{\int_{t_0}^{t}\lambda_2ds} \int_{t_0}^{t} e^{\int_{t_0}^{s}(\lambda_1-\lambda_2)d\mu}ds = -\infty \]
are satisfied; that is, \eqref{diverge3} holds. Hence, by Theorem \ref{thm:minimum1} (iii), \eqref{constoscillator} is Ulam stable on $I$, and the minimum Ulam constant is $\frac{1}{|a_0\lambda_1\lambda_2|}$. 
\end{proof}

\begin{remark}
When $a_0 = 0$, Corollary \ref{constminimum1} is completely consistent with the results given in \cite{BaiPop1,Onitsuka1}. 
\end{remark}


\section{Examples}

We now present some examples that utilize the main results of this work. In a few examples, we apply the previous theorems directly to guarantee Ulam stability of the given equation, and to find the minimal Ulam constant. In other examples, we show how the criteria of the theorems do not hold, in the cases of Ulam instability.


\begin{example}
For dimension $n=1$, consider \eqref{oscillator1} in the form of a homogeneous singular differential equation given by
\begin{equation}\label{exeq01} 
 t(1-t)x''(t) + (2-t) x'(t) + x(t) = 0, \quad t\in I=(0,1), 
\end{equation}
where $\alpha(t)=t(1-t)$, $\beta(t)=(2-t)$, and $\gamma(t)=1$ are continuous scalar functions with $\alpha(t)\ne 0$ for all $t \in I=(0,1)$. 
The associated Riccati equation \eqref{Riccati} for \eqref{exeq01} is
$$ t(1-t)\left(\rho'+\rho^2\right) + (2-t)\rho + 1 = 0, $$
which has as a solution the function
$$ \rho(t) = -\frac{1}{t}. $$
We then find that the general solution for \eqref{exeq01} with $x\left(\frac{1}{2}\right)=x_0$ and $x'\left(\frac{1}{2}\right)=x_0'$ is 
$$ x(t) = 2 x_0 + x_0' - \frac{2 x_0 + 3 x_0'}{8 t} - \frac{1}{2} \left(2 x_0 + x_0'\right) t, \quad t\in I=(0,1), $$
for arbitrary constants $x_0,x_0'\in\R$.
Using \eqref{function3} and \eqref{function4}, we calculate that both
$$ f_{3}(t) = \int_{0}^{t} e^{\int_s^t \left(-\frac{1}{\mu}\right) d\mu}ds = \frac{t}{2} $$
and
$$ f_{4}(t) = \int_{0}^{t} \frac{e^{-\int_s^t \left(-\frac{1}{\mu}+\frac{2-\mu}{\mu(1-\mu)}\right) d\mu}}{|s(1-s)|}ds 
  = \int_{0}^{t} \left(\frac{1}{s(1-s)}\right)\left(\frac{s}{1-s}\right)\left(\frac{1-t}{t}\right)ds = 1 $$
are bounded on $I=(0,1)$. Moreover, we have
$$ \lim_{t\to 0^{+}} \int_{t_0}^t \left(-\frac{1}{s}\right) ds = \lim_{t\to 0^{+}} \ln \frac{t_0}{t} = \infty $$
and
$$ \lim_{t\to 0^{+}} e^{\int_{t_0}^{t}\left(-\frac{1}{s}\right)ds} \int_{t_0}^{t} e^{-\int_{t_0}^{s}\left(-\frac{2}{\mu}+\frac{2-\mu}{\mu(1-\mu)}\right)d\mu}ds 
 = \lim_{t\to 0^{+}} \frac{t_0}{t} \int_{t_0}^{t} \frac{1-s}{1-t_0}ds = -\infty, $$
where $t_0 \in I$. 
Thus, \eqref{diverge3} holds, so that Theorems \ref{UlamS} (iii) and \ref{thm:minimum1} (iii) apply. It follows that \eqref{exeq01} is Ulam stable on $I=(0,1)$ in this case, with minimum Ulam constant
\begin{align*}
 B_{3} &:= \sup_{t \in I} \int_{0}^{t} \left(\int_{0}^{s} \frac{e^{-\int_\mu^s \left(-\frac{1}{\nu}+\frac{2-\nu}{\nu(1-\nu)}\right) d\nu}}{|\mu(1-\mu)|}d\mu\right) e^{\int_s^t \left(-\frac{1}{\mu}\right) d\mu}ds \\
&= \sup_{t \in I} \int_{0}^{t} 1 e^{\int_s^t \left(-\frac{1}{\mu}\right) d\mu}ds = \sup_{t \in I} \frac{t}{2} = \frac{1}{2}
\end{align*}
by \eqref{minimum3}. This best constant from Theorem \ref{thm:minimum1} (iii) can be verified directly. Given $\varepsilon>0$, let $x(t)=\varepsilon\left(1-\frac{t}{2}\right)$ and $\xi(t)\equiv\varepsilon$. Then $x$ is a solution of \eqref{exeq01} and $\xi$ satisfies 
$$ t(1-t)\xi''(t) + (2-t) \xi'(t) + \xi(t) = \varepsilon, \quad t\in I=(0,1), $$
with
$$ \sup_{t\in(0,1)}|\xi(t)-x(t)| = \sup_{t\in(0,1)} \frac{t\varepsilon}{2} = \frac{1}{2}\varepsilon. $$
In summary, \eqref{exeq01} is Ulam stable with minimum Ulam stability constant $B_3=\frac{1}{2}$. As can be seen by looking at the general solution, our theorems have the strength to apply to solutions that blow up at $t=0$.
\end{example}

The following four examples all deal with the Lane-Emden differential equation, either of index 0 or index 1. See also the recent paper \cite{dien}.


\begin{example}
Consider the Lane-Emden differential equation \eqref{oscillator1} given by
\begin{equation}\label{exeq02} 
 x''(t) + \frac{2}{t} x'(t) + 1 = 0, \quad t\in I=(1,\infty), 
\end{equation}
where $\alpha(t)=1$, $\beta(t)=\frac{2}{t}$, $\gamma(t)=0$, and $f(t)\equiv -1$ are continuous scalar functions with $\alpha(t)\ne 0$ for all $t \in I=(1,\infty)$. 
The associated Riccati equation \eqref{Riccati} is
$$ 1\left(\rho'+\rho^2\right) + \frac{2}{t}\rho = 0, $$
which has as a solution the function
$$ \rho(t) = -\frac{1}{t}. $$
We then find that the general solution for \eqref{exeq02} with $x(1)=x_0$ and $x'(1)=x_0'$ is 
$$ x(t) = \frac{-2+3t-t^3+6t x_0-6 x_0'+6t x_0'}{6t}, \quad t\in(1,\infty), $$
for arbitrary constants $x_0,x_0'\in\R$.
Using \eqref{function1} and \eqref{function3}, we calculate that
$$ f_{1}(t)= \int_{t}^{\infty} \frac{e^{\int_t^s \Re\left(-\frac{1}{\mu}+\frac{2}{\mu}\right) d\mu}}{|1|}ds = \int_t^{\infty}\frac{s}{t}ds = \infty $$
and
$$ f_{3}(t)= \int_{1}^{t} e^{\int_s^t \Re(-\frac{1}{\mu}) d\mu}ds = \frac{t}{2}-\frac{1}{2t} $$
is unbounded on $I=(1,\infty)$, so that Theorem \ref{UlamS} does not apply. Indeed, given an arbitrary $\varepsilon>0$,
$$ \xi''(t) + \frac{2}{t} \xi'(t) + 1 = \varepsilon, \quad t\in I=(1,\infty) $$
has a solution $\xi(t)=\frac{(\varepsilon-1)t^2}{6}$, and thus
$$ \sup_{t\in I}|\xi(t)-x(t)| = \sup_{t\in I}\left|\frac{2 + 6 x_0' - 3 t (1 + 2 x_0 + 2 x_0') + t^3 \varepsilon}{6t}\right|=\infty $$
for any choice of $x_0,x_0'\in\R$, so that in fact \eqref{exeq02} is not Ulam stable on $(1,\infty)$. 
\end{example}


\begin{example}\label{example3}
In this example we modify the interval $I$ for the Lane-Emden equation \eqref{exeq02} by considering
\begin{equation}\label{exeq03} 
 x''(t) + \frac{2}{t} x'(t) + 1 = 0, \quad t\in I=(0,\sigma), \quad 0<\sigma<\infty. 
\end{equation}
It is easy to verify that $f_3$ in \eqref{function3} and $f_4$ in \eqref{function4} have expression
$$ f_3(t)=\frac{t}{2}=f_4(t), \quad t\in(0,\sigma), $$
as $\rho(t)=-\frac{1}{t}$, $\alpha(t)=1$, and $\beta(t)=-2\rho(t)$. Moreover, we have
$$ \lim_{t\to 0^{+}} \int_{t_0}^t \left(-\frac{1}{s}\right) ds = \lim_{t\to 0^{+}} \ln \frac{t_0}{t} = \infty $$
and
$$ \lim_{t\to 0^{+}} e^{\int_{t_0}^{t}\left(-\frac{1}{s}\right)ds} \int_{t_0}^{t} e^{-\int_{t_0}^{s}\left(-\frac{2}{\mu}+\frac{2}{\mu}\right)d\mu}ds 
 = \lim_{t\to 0^{+}} \frac{t_0}{t} (t-t_0) = -\infty, $$
where $t_0 \in I$. Then \eqref{diverge3} holds, so that \eqref{exeq03} is Ulam stable on $I=(0,\sigma)$ for any $\sigma\in(0,\infty)$, with minimum Ulam constant
\begin{align*}
 B_{3} &:= \sup_{t \in I} \int_{0}^{t} \left(\int_{0}^{s} e^{-\int_\mu^s \left(\frac{1}{\nu}\right) d\nu}d\mu\right) e^{\int_s^t \left(-\frac{1}{\mu}\right) d\mu}ds \\
&= \sup_{t \in (0,\sigma)} \int_{0}^{t} \left(\frac{s}{2}\right)\left(\frac{s}{t}\right)ds = \sup_{t \in I} \frac{t^2}{6} = \frac{\sigma^2}{6}
\end{align*}
by \eqref{minimum3}, after employing Theorems \ref{UlamS} (iii) and \ref{thm:minimum1} (iii). In summary, \eqref{exeq03} is Ulam stable on $(0,\sigma)$ for any $\sigma\in(0,\infty)$, with minimum Ulam stability constant $B_3=\frac{\sigma^2}{6}$. 
\end{example}


\begin{example}
Consider \eqref{oscillator1} in the form of the Lane-Emden differential equation given by
\begin{equation}\label{exeq04} 
 x''(t) + \frac{2}{t} x'(t) + x(t) = 0, \quad t\in I=(0,\sigma), \quad \sigma\in\left(0,\frac{\pi}{2}\right), 
\end{equation}
where $\alpha(t)=1$, $\beta(t)=\frac{2}{t}$, $\gamma(t)=1$, and $f(t)\equiv 0$ are continuous scalar functions with $\alpha(t)\ne 0$ for all 
$t \in I=(0,\sigma)$, where $\sigma\in\left(0,\frac{\pi}{2}\right)$. 
The associated Riccati equation \eqref{Riccati} is
$$ 1\left(\rho'+\rho^2\right) + \frac{2}{t}\rho + 1 = 0, $$
which has as a solution the function
$$ \rho(t) = -\tan(t)-\frac{1}{t}, \quad t \in(0,\sigma).  $$
We then find that the general solution for \eqref{exeq04} is 
$$ x(t) = \frac{c_1\cos t}{t} + \frac{c_2\sin t}{t}, \quad t\in(0,\sigma), $$
for arbitrary constants $c_1,c_2\in\R$.
Using \eqref{function3} and \eqref{function4}, we calculate that both
$$ f_{3}(t) = \int_{0}^{t} e^{\int_s^t (-\tan(\mu)-\frac{1}{\mu}) d\mu}ds = \int_{0}^{t} \frac{s\cos(t)}{t\cos(s)}ds $$
and
$$ f_{4}(t) = \int_{0}^{t} e^{-\int_s^t \left(-\tan(\mu)+\frac{1}{\mu}\right) d\mu}ds = \tan(t)+\frac{\cos(t)-1}{t\cos(t)} $$
are bounded on $I=(0,\sigma)$, because
$$ \lim_{t\to 0^{+}} \int_{0}^{t} \frac{s\cos(t)}{t\cos(s)}ds = 0 $$
holds. Moreover, we have
$$ \lim_{t\to 0^{+}} \int_{t_0}^t \left(-\tan(s)-\frac{1}{s}\right) ds = \lim_{t\to 0^{+}} \ln \frac{t_0\cos(t)}{t\cos(t_0)} = \infty $$
and
\begin{align*}
 \lim_{t\to 0^{+}} e^{\int_{t_0}^{t}\left(-\tan(s)-\frac{1}{s}\right)ds} \int_{t_0}^{t} e^{-\int_{t_0}^{s}\left[2\left(-\tan(\mu)-\frac{1}{\mu}\right)+\frac{2}{\mu}\right]d\mu}ds 
 &= \frac{t_0\cos(t)}{t\cos(t_0)} \int_{t_0}^{t} \frac{\cos^2(t_0)}{\cos^2(s)}ds\\
 &= \frac{\sin(t-t_0)}{t} 
 = -\infty,
\end{align*}
where $t_0 \in I$. Then \eqref{diverge3} holds, so that Theorems \ref{UlamS} (iii) and \ref{thm:minimum1} (iii) apply. It follows that \eqref{exeq04} is Ulam stable on $I=(0,\sigma)$ in this case, with minimum Ulam constant
\begin{align*}
 B_{3} &:= \sup_{t \in I} \int_{0}^{t} \left(\int_{0}^{s} e^{-\int_\mu^s \left(-\tan(\nu)-\frac{1}{\nu}+\frac{2}{\nu}\right) d\nu}d\mu\right) e^{\int_s^t \left(-\tan(\mu)-\frac{1}{\mu}\right) d\mu}ds \\
  &= \sup_{t \in I} \int_{0}^{t} \left(\tan(s)+\frac{\cos(s)-1}{s\cos(s)}\right)\left(\frac{s\cos(t)}{t\cos(s)}\right)ds\\
  &= \sup_{t \in I} \left(\frac{\cos(t)}{t}\right)\left(\frac{t-\sin(t)}{\cos(t)}\right) = 1-\frac{\sin(\sigma)}{\sigma}
\end{align*}
by \eqref{minimum3}. In summary, \eqref{exeq04} is Ulam stable on $(0,\sigma)$ for any $\sigma\in\left(0,\frac{\pi}{2}\right)$, with minimum Ulam stability constant $B_3=1-\frac{\sin(\sigma)}{\sigma}$. 
\end{example}


\begin{example}
In this example we modify the interval $I$ for the Lane-Emden equation \eqref{exeq04} by considering
\begin{equation}\label{exeq05} 
 x''(t) + \frac{2}{t} x'(t) + x(t) = 0, \quad t\in I=(\tau,\infty), \quad 0\le\tau<\infty. 
\end{equation}
Once again, the corresponding Riccati equation has solution $\rho(t)=-\tan(t)-\frac{1}{t}$, which does not exist infinitely often on $I=(\tau,\infty)$ due to the tangent function, so Theorems \ref{UlamS} and \ref{thm:minimum1} cannot be applied. Given $\varepsilon>0$, consider 
$$\xi(t)=\frac{\varepsilon}{8t}\left(\cos(t)+2t\sin(t)-2t^2\cos(t)\right), \quad t\in(\tau,\infty).$$
Note that
$$ \sup_{t \in I} \left|\xi''(t)+\frac{2}{t}\xi'(t)+\xi(t)\right| = \sup_{t \in I} \left|\varepsilon\sin(t)\right| = \varepsilon. $$
Since 
$$ x(t) = \frac{c_1\cos t}{t} + \frac{c_2\sin t}{t}, \quad t\in(\tau,\infty) $$
is the general solution for \eqref{exeq05}, we have
$$ \sup_{t \in I} |\xi(t)-x(t)| =\infty, $$
making \eqref{exeq05} unstable in the Ulam sense on $I=(\tau,\infty)$ for any $\tau\in[0,\infty)$.
\end{example}


\begin{example}
In this example, we consider an extension of Example \ref{example3} to include more general power functions. We consider the second-order linear differential equation
\begin{equation}\label{exeq06} 
 t^{1-a}x''(t) + bt^{-a}x'(t) + (b-2)t^{-1-a}x(t) + t^{b-2} = 0, \quad t\in I=(0,\sigma), \quad 0<\sigma<\infty, 
\end{equation}
where $a$ and $b$ are real-valued constants with
$$ 1-a < b \le 2. $$
If $a=1$ and $b=2$, then this equation reduces to \eqref{exeq03}. It is easy to verify that $\rho(t)=-\frac{1}{t}$ is a solution of the associated Riccati equation
$$ t^{1-a}\left(\rho'+\rho^2\right) + bt^{-a}\rho + (b-2)t^{-1-a} = 0, $$
and $f_3$ in \eqref{function3} and $f_4$ in \eqref{function4} have expression
$$ f_3(t)=\frac{t}{2}, \quad f_4(t)=\frac{t^{a}}{a+b-1}, \quad t\in(0,\sigma), $$
as $\alpha(t)=t^{1-a}$, and $\beta(t)=bt^{-a}$. Moreover, we have
$$ \lim_{t\to 0^{+}} \int_{t_0}^t \left(-\frac{1}{s}\right) ds = \lim_{t\to 0^{+}} \ln \frac{t_0}{t} = \infty $$
and
$$ \lim_{t\to 0^{+}} e^{\int_{t_0}^{t}\left(-\frac{1}{s}\right)ds} \int_{t_0}^{t} e^{-\int_{t_0}^{s}\left(-\frac{2}{\mu}+\frac{b}{\mu}\right)d\mu}ds 
 = \lim_{t\to 0^{+}} \frac{t_0}{t} \int_{t_0}^t \left(\frac{s}{t_0}\right)^{2-b} ds = \lim_{t\to 0^{+}} \frac{t_0^{b-1}}{3-b} \left(t^{2-b}-\frac{t_0^{3-b}}{t}\right) = -\infty, $$
where $t_0 \in I$. Then \eqref{diverge3} holds, so that \eqref{exeq06} is Ulam stable on $I=(0,\sigma)$ for any $\sigma\in(0,\infty)$, with minimum Ulam constant
\begin{align*}
 B_{3} &:= \sup_{t \in I} \int_{0}^{t} \left(\int_{0}^{s} \frac{e^{-\int_\mu^s \left(-\frac{1}{\nu}+\frac{b}{\nu}\right) d\nu}}{|s^{1-a}|}d\mu\right) e^{\int_s^t \left(-\frac{1}{\mu}\right) d\mu}ds \\
&= \sup_{t \in (0,\sigma)} \int_{0}^{t} \left(\frac{s^{a}}{a+b-1}\right)\left(\frac{s}{t}\right)ds = \sup_{t \in I} \frac{t^{a+1}}{(a+2)(a+b-1)} = \frac{\sigma^{a+1}}{(a+2)(a+b-1)}
\end{align*}
by \eqref{minimum3}, after employing Theorems \ref{UlamS} (iii) and \ref{thm:minimum1} (iii). In summary, \eqref{exeq06} is Ulam stable on $(0,\sigma)$ for any $\sigma\in(0,\infty)$, with minimum Ulam stability constant $B_3=\frac{\sigma^{a+1}}{(a+2)(a+b-1)}$. 
\end{example}


\section{Conclusions}
This work investigated the Ulam stability of second-order linear differential vector equations with variable coefficients. Sufficient conditions for Ulam stability and explicit Ulam stability constants are given. In particular, if restricted to real-valued coefficients, the best Ulam constants are derived. To the best of the authors' knowledge, no best Ulam constants are known so far for second-order non-autonomous equations other than periodic systems. Therefore, this is the first study to derive best Ulam constants for second-order non-periodic non-autonomous linear differential equations. Various non-trivial examples, mainly Lane-Emden differential equations, were provided to illustrate the results obtained. Some examples show that Ulam stability can be guaranteed even for solutions that blow up in finite time, and that best Ulam constants can be derived. Note that it is emphasized here that it is also the first time that a Ulam stability analysis of blow--up solutions of the second-order equation has been presented. In addition, examples of instability are also presented.


\section*{Acknowledgments}
M. O. was supported by the Japan Society for the Promotion of Science (JSPS) KAKENHI (grant number JP20K03668).

\end{document}